\documentclass[12pt]{amsart}
\usepackage{amscd,amssymb,epsfig}
\usepackage{graphics}

\input epsf
\newtheorem{theorem}{Theorem}[section]
\newtheorem{lemma}[theorem]{Lemma}
\newtheorem{proposition}[theorem]{Proposition}
\newtheorem{corollary}[theorem]{Corollary}
\theoremstyle{definition}
\newtheorem{definition}[theorem]{Definition}

\newtheorem{example}[theorem]{Example}

\theoremstyle{remark}

\numberwithin{figure}{section}
\numberwithin{table}{section}
\def\varph{{\varphi}}
\def\ra{{\rightarrow}}
\def\lra{{\longrightarrow}}
\def\om{{\omega}}

\def\La{{\Lambda}}
\def\ga{{\gamma}}
\def\Ga{{\Gamma}}
\def\vGa{{\varGamma}}

\def\bZ{{\mathbb Z}}
\def\bQ{{\mathbb Q}}

\def\Ig{{\mathcal I}_g}
\def\Isg{{\mathcal I}_{g,*}}
\def\Msg{{\mathcal M}_{g,*}}
\def\M1g{{\mathcal M}_{g,1}}
\def\I1g{{\mathcal I}_{g,1}}
\newcommand\Hom{\operatorname{Hom}}

\newcommand\sgn{\operatorname{sgn}}

\catcode`\@=\active 

\begin{document}

\title[Torelli groups, Johnson homomorphisms, and new cycles]
{Torelli groups, extended Johnson homomorphisms,
and new cycles on the moduli space of curves}
\author{S.~Morita}
\address{Department of Mathematical Sciences\\
University of Tokyo \\Komaba, Tokyo 153-8914\\
Japan}
\email{morita{\char'100}ms.u-tokyo.ac.jp}
\author{R. C. Penner}
\address{Departments of Mathematics and Physics/Astronomy\\
University of Southern California\\
Los Angeles, CA 90089\\
USA}
\email{rpenner{\char'100}math.usc.edu}
\thanks{The first author is partially supported by JSPS Grant 16204005}

\begin{abstract}
Infinite presentations are given for all of the higher Torelli groups
of once-punctured surfaces. In the case of the classical Torelli
group, a {\it finite presentation} of
the corresponding {\it groupoid} is also given, and finite 
presentations of the classical Torelli groups acting trivially
on homology modulo $N$ are derived for all $N$.
Furthermore, the first Johnson homomorphism, which is defined from
the classical Torelli group to the third exterior power of the
homology of the surface, is shown to {\it lift} to an
explicit {\it canonical} $1$-{\it cocycle} of the Teichm\"uller space.
The main tool for these results is the known mapping class group
invariant ideal cell decomposition of the Teichm\"uller space.

This new $1$-cocycle is mapping class group equivariant, so
various contractions of its powers yield
various combinatorial (co)cycles of the moduli space
of curves, which
are also new.
Our combinatorial construction can be related to
former works of Kawazumi and the first-named author with the consequence
that the algebra generated by the cohomology classes
represented by the new cocycles is precisely the tautological
algebra of the moduli space.

There is finally a discussion of
prospects for similarly finding cocycle lifts of
the higher Johnson homomorphisms.

\end{abstract}

\maketitle

\section{Introduction}

Fix a smooth, closed, oriented surface $F$ of genus $g\geq 1$
and also fix a base point $*\in F$.
We consider the mapping class group
$$
\Msg =\pi_0\, \mathrm{Diff}_+ (F,*)
$$
of $F$ relative to the base point $*\in F$.
Let $\pi=\pi_1(F,*)$
denote the fundamental group with abelianization
the first integral homology group $H=H_1(F;\bZ)$.
The group $\pi$ has lower central series defined recursively by
$\Ga_0=\pi$ and $\Ga_{k+1}=[\pi,\Ga_k]$ for $k\geq 0$, where the bracket
of two groups denotes their commutator group, and the k-th nilpotent
quotient is defined as
$N_k=\pi/\Ga_k$. Thus, we have the exact sequence
\begin{equation}
0\lra \Ga_k/\Ga_{k+1}\lra N_{k+1}\lra N_k\lra 1
\label{eq:nk}
\end{equation}
which is a central extension, and by a result of Labute \cite{Labute},
the quotient $\Ga_k/\Ga_{k+1}$
is identified with the degree $k+1$ part of the free Lie algebra
$\mathcal{L}(H)$ of $H=N_1$ divided by the ideal generated by the symplectic
class $\om_0\in \mathcal{L}_2(H)=\Lambda^2 H$ in degree two.

The inverse limit $N_\infty$ of the tower
$\cdots \ra N_{k+1}\ra N_k\ra\cdots\ra N_1\ra 1$
is the pronilpotent completion of $\pi$. In the spirit
of rational homotopy theory \cite{Sullivan77} or
Mal'cev theory \cite{Malcev}, $N_\infty$ or its
rational model should be thought of as a kind of
approximation to $\pi$ itself.

A mapping class $\varph\in\Msg$ induces actions $\varph_k:N_k\ra N_k$,
for each $k\geq 1$ and the kernel of the homomorphism
$\Msg\ra\mathrm{Aut}(N_k)$ is the {\it k-th Torelli group} to be
denoted $\Isg(k)\subset\Msg$.
In particular, for $k=1$, we find the classical Torelli group
$0\ra\Isg\ra\Msg\ra\mathrm{Sp}(2g,\bZ)\ra 1$,
where $\mathrm{Sp}(2g,\bZ)$ denotes the Siegel modular group.

If $\varph\in\Msg$ acts trivially on $N_k$ and $\gamma\in N_{k+1}$,
then $\varph(\gamma)\gamma^{-1}\in \mathrm{Ker}(N_{k+1}\ra N_k)$, so by
exactness of
\eqref{eq:nk} together with a few arguments, there are mappings
$$
\tau_k: \Isg(k)\lra\mathrm{Hom}(N_{k+1}, \Ga_k/\Ga_{k+1}),
$$
and these are called the Johnson homomorphisms
introduced in \cite{Johnson80}\cite{Johnson83a}
(see also prior works of Andreadakis
\cite{Andreadakis} and Sullivan \cite{Sullivan75}).
See \cite{Johnson83a} for a survey of the Torelli groups
and \cite{Morita93a}\cite{Morita99}\cite{Morita05} for further results.

In particular for $k=1$, Johnson \cite{Johnson83}
has given a finite generating set for each of the Torelli groups $\Isg,\Ig$
($\Ig$ denotes the Torelli group for the closed surface $F$)
when the genus is at least three. He has also calculated the kernel
of $\tau_1$ in \cite{Johnson85} and described the abelianization
of $\Ig$ in \cite{Johnson85b}.
In genus two, Mess \cite{Mess} has shown that $\mathcal{I}_2$ is
an infinitely generated free group.  Biss and Farb \cite{BissFarb} have also
shown that the second Torelli groups
$\mathcal{I}_g(2)$ and $\mathcal{I}_{g,*}(2)$,
which are usually denoted by $\mathcal{K}_g$
and $\mathcal{K}_{g,*}$,
are not finitely generated for any $g\geq 2$.

In this paper, we give the first presentations of all the
higher Torelli groups $\Isg(k)$ for $\infty > k\geq 1$ (see $\S5$),
but they have
infinitely many generators and infinitely many relations.
As part of our approach, we shall also study corresponding
``Torelli groupoids" (whose definition will be given in $\S2$), which also
admit infinite presentations as shown in $\S5$.
In the classical case $k=1$ using Johnson's finite generating
set, we furthermore give a finite presentation of the classical
Torelli {\it groupoids} in $\S2$.  Finite presentability of the corresponding
classical Torelli groups is a long-standing problem.

The key tool for these investigations is the mapping class group
invariant {\it ideal} cell decomposition \cite{Strebel}\cite{Hubbard-Masur}
\cite{Harer85}\cite{Harer86} (in the conformal context) and
\cite{Penner87}\cite{Penner88}\cite{Penner92}\cite{BE}\cite{Penner96} (in 
the hyperbolic context)
of the Teichm\"uller space
$\mathcal{T}_{g,*}$
of $F$ relative to the base point $*$, where $*$ is often regarded as a 
puncture geometrically
(see the next section).
By equivariance, there is an induced ideal cell decomposition of the 
corresponding
moduli space $\mathbf{M}_{g,*}$,
and this cell decomposition can be used to give a finite
essentially combinatorial presentation of the mapping class group itself.
The idea is very simple: any path can be put into general position
with respect to the codimension one faces of the cell decomposition,
so generators are given by crossing these faces; any homotopy of paths
can be put into general position with respect to the codimension two
faces of the cell decomposition, and an analysis of the links
of these faces then provides a complete set of relations.
(It is actually a bit more complicated as
${\mathbf M}_{g,*}$ is an orbifold, cf.  \cite{Igusa2}.)

Letting $\mathbf{T}_k=\mathcal{T}_{g,*}/\Isg(k)$ denote the
{\it k-th Torelli space}, there is thus a tower
\begin{equation}
\mathcal{T}_{g,*}\lra\cdots\lra
\mathbf{T}_{k+1}\lra \mathbf{T}_k\lra \cdots \mathbf{T}_1\lra
\mathbf{M}_{g,*}
\label{eq:spaces}
\end{equation}
of spaces and covering maps.
The ideal cell decomposition of Teichm\"uller space
$\mathcal{T}_{g,*}$ descends to ideal cell decompositions on each
Torelli space $\mathbf{T}_k$, and this gives a useful combinatorial
model of the tower \eqref{eq:spaces}.

The cell decomposition also leads to a new combinatorial expression
for the classical Johnson homomorphism
$\tau_1:\Isg\ra\La^3 H$. We shall find that
$\tau_1$ descends from an explicit $1$-cocycle
$j\in Z^1(\hat{\mathcal G}_{T};\La^3 H)$, where $\hat{\mathcal G}_T$
denotes the Poincar\'e dual to the ideal cell decomposition of
$\mathcal{T}_{g,*}$. More precisely, $\hat{\mathcal G}_T$ is a
(genuine) cell complex of dimension $4g-3$ each of whose $k$-cell
is Poincar\'e dual to a $(6g-4-k)$-cell of $\mathcal{T}_{g,*}$.
Furthermore $\hat{\mathcal G}_T$ sits inside $\mathcal{T}_{g,*}$
and there exists an $\Msg$-equivariant deformation retraction of the
whole space onto this subspace (see \cite{Harer86} for more details).

The above $1$-cocycle $j$ is $\Msg$-equivariant
and we can show that $\tau_1$ is the restriction of a crossed
homomorphism from $\Msg$ to $\La^3 H$. This gives a new proof that
$\tau_1$ so extends, as first shown in \cite{Morita93}.

The group theoretical method of \cite{Morita96} can be adapted
to this combinatorial setting and,
as a consequence, we obtain explicit (co)cycles of the moduli
space $\mathbf{M}_{g,*}$.

This paper is organized as follows. The mapping class group
invariant ideal cell decomposition of the Teichm\"uller space
$\mathcal{T}_{g,*}$ is recalled in $\S2$,
and the infinite presentation of $\Isg$, the finite presentation
of the corresponding classical Torelli {groupoid}, and the finite
presentations of the classical Torelli groups modulo $N$ are given.
In $\S3$, we define and confirm
the putative combinatorial expression for the Johnson
homomorphism $\tau_1$ as a $1$-cocycle $j\in Z^1(\hat{\mathcal G}_T;\La^3 H)$.
In $\S4$, we construct various (co)cycles of the moduli space by making use
of $j$ and describe the cohomology classes represented by them.
Finally, in $\S5$ we give the infinite presentations of
higher Torelli groups in analogy to the treatment
of $\Isg$ in $\S2$ and give a combinatorial description but {\it not a 
cocycle} for the higher
Johnson homomorphisms. Concluding remarks and a discussion of further 
prospects towards analogous
cocycle lifts of the higher Johnson homomorphisms are given in
$\S6$.
\par
\vspace{1cm}

\noindent
{\bf Acknowledgements}\quad
It is a pleasure for the second-named author to thank the University of 
Tokyo for hospitality
during late 2004 when this work was initiated, the kind support
of the Fields Institute where part of this work was completed, and
Bob Guralnick for numerous useful discussions.  Finally, let us thank
Matthew Day for pointing out Corollary~\ref{levelN} and Nate Broaddus
for pointing out Corollary~\ref{secondtor}.

\vskip .2in

\section{Presentations of the classical Torelli group and groupoid}

Fix a smooth closed oriented surface $F$ of genus $g\geq 1$ with a
base point (or ``puncture'') $*\in F$,
and let
${\mathcal T}_{g,*}={\mathcal T}(F')$ denote the {\it Teichm\"uller space} of
complete finite-area metrics of constant
curvature -1 on $F'=F-\{ *\}$ modulo push-forward by diffeomorphisms
isotopic to the
identity.  Thus, the mapping class group $\Msg$ acts properly
discontinuously with fixed points
on $\mathcal{T}_{g,*}$, and the quotient orbifold is {\it (Riemann's)
moduli space}
$\mathbf{M}_{g,*}={\mathbf M}(F,*)$.

The action of
$\Msg$ on the first integral relative homology group
$$
H=H_1(F,\{ *\} )\cong H^1(F,\{ *\})\cong H_1(F-\{ *\})
$$
gives a natural map
$\Msg\to \mathrm{Sp}(2g,\bZ)$ onto the Siegel modular group,
and the kernel of this map is the {\it (classical) Torelli group}
$\Isg \subset \Msg$, with quotient manifold
{\it Torelli space} $\mathbf{T}_{g,*}=\mathcal{T}(F')/\Isg$.

According to \cite{Penner88}, the Teichm\"uller space
$\mathcal{T}_{g,*}$
comes equipped with a canonical $\Msg$-invariant
ideal cell decomposition ${\mathcal G}_T$, where the cells in ${\mathcal G}_T$
are in one-to-one correspondence with isotopy classes of fat graph spines
embedded in $F'$
and where the face relation is induced by the collapse of edges with
distinct endpoints.  This
``fat graph complex'' ${\mathcal G}_T$ descends to an ideal cell decomposition
${\mathcal G}_M={\mathcal G}_T/\Msg$ of
$\mathbf{M}_{g,*}$ and indeed also descends to an ideal cell decomposition
${\mathcal G}_I={\mathcal G}_T/\Isg$ of the
Torelli space $\mathbf{T}_{g,*}$ as well.
Thus, the tower $\mathcal{T}_{g,*}\to \mathbf{T}_{g,*}\to \mathbf{M}_{g,*}$
of spaces admits the combinatorial model
${\mathcal G}_T\to {\mathcal G}_I\to{\mathcal G}_M$ of ideal cell complexes.
Notice that $\mathbf{T}_{g,*}\to \mathbf{M}_{g,*}$ is an
infinite-sheeted cover branched over the non-manifold points of
$\mathbf{M}_{g,*}$ with deck group
$\mathrm{Sp}(2g,\bZ)$.

Poincar\'e dual to the fat graph spine $G$ of the punctured surface $F'$ is
a collection of arcs
$A _G$ connecting $*$ to itself which decompose $F$ into polygons.  It will
sometimes be
useful to employ arc families rather than fat graphs.  The dual of the
collapse of
an edge of $G$ is the removal of its dual arc from $A _G$.

It is also convenient to take the Poincar\'e dual complexes
$\hat{\mathcal G}_T\to \hat{\mathcal G}_I\to\hat{\mathcal G}_M$,
where there is one $k$-dimensional cell of $\hat{\mathcal G}_X$ for
each ($6g-4-k$)-dimensional cell of ${\mathcal G}_X$, for $X=M,I,T$.  In
particular, the
zero-dimensional cells of $\hat{\mathcal G}_T$ are in one-to-one
correspondence with
isotopy classes of trivalent fat graph spines of $F'$, and the
one-dimensional cells of
$\hat{\mathcal G}_T$ are in one-to-one correspondence with classes of fat
graphs which are
trivalent except for one four-valent vertex, i.e.,
the oriented one-dimensional cells are in one-to-one correspondence
with pairs consisting of a trivalent fat graph together with a distinguished
edge of it with distinct endpoints.
There is furthermore
a two-dimensional cell in $\hat{\mathcal G}_T$ for each class of fat graph
which is trivalent
except for either one five-valent vertex or two four-valent vertices.

The (orbifold) fundamental group of $\hat{\mathcal G}_X$ is  respectively
$\Msg$, $\Isg$,
trivial,  for $X=M,I,T$; we shall also require the fundamental path
groupoids of these spaces (i.e., objects are vertices of $\hat{\mathcal
G}_X$ and morphisms are
homotopy classes of paths) to be denoted $\Gamma _X$.  In particular, 
$\Isg$ is
identified with
the subgroup of $\Gamma_I$ with both endpoints at some fixed vertex.

Since the action of the Torelli group $\Isg$ on the Teichm\"uller space
is known to be {\it free}, we can use the well known procedure to obtain
a presentation of $\Isg$:
generators are given by ``Whitehead moves''
(see Figure~1) along edges in $\Isg$-orbits of
fat graphs and relations are given by $\Isg$-orbits of codimension-two cells.

A {\it homology marking} on a trivalent fat graph $G$ is the assignment of
a class
$a\in H$ to each oriented edge $a$ of $G$ (thus abusing notation slightly) so
that:
\begin{enumerate}
\item
if $a,b$ are the two
orientations on a common edge, then $a=-b$;
\item
if three distinct oriented edges
$a,b,c$ point towards a common vertex of $G$, then $a+b+c=0$;
\item
$\{ a\in H:a~{\rm is~an~oriented~edge~of}~G\}$ spans $H$, that is, the
assignment has full rank.
\end{enumerate}

For each class of embedding of trivalent fat graph spine $G$ in $F'$ (i.e.,
for each
name of a $(6g-4)$-cell in ${\mathcal G}_T$),  there is a canonical
homology marking on $G$, namely,
Poincar\'e dual to $G$ in $F'$ is an arc family $\alpha$, where there is
one arc in
$\alpha$ connecting
the puncture $*$ to itself dual to each edge of $G$.  Given an orientation
$a$ on this edge of
$G$, there is a unique orientation to its dual arc so that the pair of
orientations (in some
fixed order, say, first the dual arc and then the edge of $G$) are
compatible with the
orientation of
$F$ itself.  We may assign to the edge $a$ the homology class of the dual
arc so oriented;
property 1) follows by definition, property 2) since the three dual arcs
bound a triangle in
$F$ rel $\{ *\} $, and property 3) since $G$ is a spine of $F'$.


\vskip .2in

\centerline{\epsffile{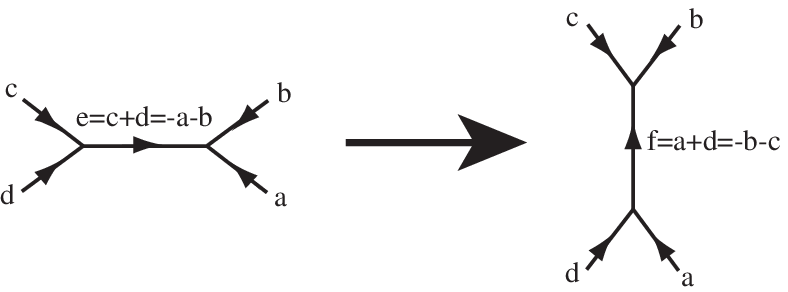}}

\vskip .2in

\centerline{{\bf Figure 1}~Whitehead move on homology marking.}

\vskip .2in

A Whitehead move acts on a homology marking in the natural way as
illustrated in
Figure~1.  A sequence of Whitehead moves corresponding to a loop in
$\mathbf{M}_{g,*}$ begins and ends with
fat graphs which are related by a mapping class $\varph:(F,*)\to (F,*)$;
furthermore, there is the
combinatorial identification of edges $e$ with $f$ leaving the other edges
unchanged, as in
Figure~1, and
$\varph$ represents an element of $\Isg$ if and only if under this
combinatorial identification of edges, the homology marking is furthermore
invariant under the
sequence of Whitehead moves.

Let $W_{(G,e)}$, or just $W_e$ if $G$ is fixed or understood, denote the 
Whitehead move along
the unoriented edge $e$ of $G$, where we always assume that $e$ has
distinct endpoints.

\begin{theorem}
Fix an isotopy class of trivalent fat graph spine $G$ in $F'$.
Then the Torelli group
$\Isg$ is generated by
sequences of Whitehead moves starting from
$G$ and ending at
$\varph (G)$ for some $\varph\in \Msg$, where the sequence of Whitehead
moves leaves invariant the
homology marking.  Furthermore, a complete set of relations in $\Isg$ is
generated by:
\begin{enumerate}
\item
{\rm involutivity relation}:
$W_e W_f=1=W_f W_e$ in the notation of Figure~1;
\item
{\rm commutavity relation}: if $e,f$ share no vertices, then $W_e$ and $W_f$
commute;
\item
{\rm pentagon relation}:
$W_f W_{g_1} W_{f_2} W_{g_3} W_{f_4} =1$ in the circumstance and notation of
Figure~2.
\end{enumerate}
\label{th:presentation}
\end{theorem}


\vskip .2in

\centerline{\epsffile{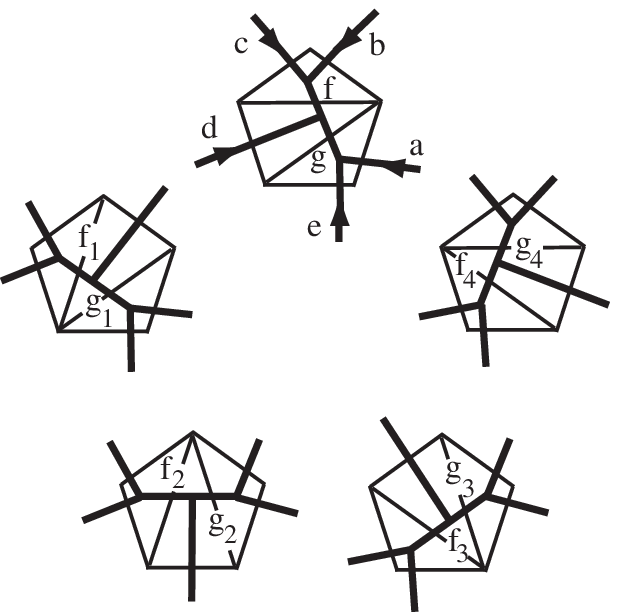}}

\vskip .2in

\centerline{{\bf Figure 2}~The pentagon relation.}

\vskip .2in

\begin{proof}
By definition, cells in ${\mathcal G}_I$ are indexed by $\Isg$-orbits of
homology marked fat graphs, and since the natural map $\Msg\to 
Sp(2g,{\mathbb Z})$ is a
surjection, every $\Isg$-orbit of homology marked trivalent fat graph 
actually arises.
Given
$\varph\in\Isg$, we may choose a point  in $\xi\in \mathbf{T}_{g,*}$ and a
path connecting $\xi$ to $\varph (\xi )$.  This path may be perturbed so as
to meet the
codimension-one cells transversely and hence be expressed as a collection
of Whitehead moves.
Furthermore,
$\varph\in\Isg$ precisely means that the sequence of Whitehead moves leaves
invariant some
(hence any) homology marking.  Thus, the set of homology marking-preserving
sequences of
Whitehead moves beginning and ending on the same $\Msg$-orbit
is not only a generating set for
$\Isg$, it is in fact an enumeration (with repeats) of $\Isg$.

The relations likewise follow from general position of
codimension-two cells of ${\mathcal G}_I$ with path homotopies in
$\mathbf{T}_{g,*}$.  Indeed, involutivity
corresponds to consecutively crossing the same codimension-one cell, and
the pentagon and
commutativity relations correspond to the two types of codimension-two
cells in ${\mathcal G}_T$
or ${\mathcal G}_I$.
\end{proof}

Thus, the Torelli group admits a presentation akin to that of the mapping
class group.  There are also analogous infinite presentations of all the
higher Torelli groups
which will be given in $\S$5.

On the other hand special to the
case of the classical Torelli groups, there is Johnson's finite generating
set of $\Isg$,
which we can combine with the combinatorics of the cell decomposition to give a
{\sl finite presentation} of the ``Torelli groupoid'' $\Gamma _I$.

\begin{theorem}\label{groupoid}
Suppose that the genus $g$  of $F$ is at least $3$.  Choose a
cellular fundamental domain
$D\subseteq \mathcal{T}_{g,*}$ for the action of $\Msg$ on $\mathcal{T}_{g,*}$,
and choose a top-dimensional cell in $D$ with
corresponding class of fat graph spine $G$ in $F'$.  Then
$\Gamma _I$ is generated by:
\begin{enumerate}
\item
for each codimension-one face in $D$, there is a generator given by the
corresponding Whitehead move;
\item
for each of Johnson's generators of $\Isg$, choose a representing
sequence of Whitehead moves starting from $G$;
\item
for each of the five generators of $\mathrm{Sp}(2g,\bZ)$, choose a
representing sequence of Whitehead moves starting from $G$.
\end{enumerate}
Furthermore, a complete set of relations in $\Isg$ is given by
involutivity and the commutativity or pentagon relations corresponding to
the codimension-two
faces in $D$.
\label{th:2}
\end{theorem}

\begin{proof}
Let $B$ denote the groupoid generated by the Whitehead moves or sequences
of Whitehead moves in
1-3), so necessarily $B<\Gamma _I$.  By general position as before, it
suffices both for
generators and relations to prove that the
$B$-orbits of cells in $D$ cover $\mathcal{T}_{g,*}$.  For if this is so in
codimension-one, then the
generators of $B$ include the generators of $\Gamma _I$, so $B=\Gamma _I$,
and if this is so
in codimension-two, then all codimension-two cells of ${\mathcal G}_I$
(i.e., all the relations
in $\Gamma _I$) arise in the
$B$-orbits of the codimension-two cells in $D$.  Thus, all relations in
$\Gamma _I$ lie in
the $B$-normal closure of the asserted relations.

The proof that $B(D)=\mathcal{T}_{g,*}$ is proved using standard covering
space theory from the
facts that
$\mathbf{T}_{g,*}\to \mathbf{M}_{g,*}$ is a covering map,
the deck group $\mathrm{Sp}(2g,\bZ)$
acts transitively on the fibers of
$\mathbf{T}_{g,*}\to \mathbf{M}_{g,*}$,
and generators for the fundamental group of $\mathbf{T}_{g,*}$ are supplied
by 2).

In more detail and with more notation, $B$ is generated by: 1) the
Whitehead moves
$W_{(G_0,e_0)}$ for each trivalent fat graph $G_0$ with cell in $D$ and
edge $e_0$ of $G_0$
with distinct endpoints, 2) Johnson's generators $u_1,\ldots ,u_L$ for
$\Isg$, and 3)
generators $v_1,\cdots ,v_5$ for $\mathrm{Sp}(2g,\bZ)$.

Suppose that $\gamma\in \Gamma _I$ is a sequence of Whitehead moves
starting from
the fat graph $G_0$ and ending at the fat graph $G_1$.  Under the natural
identification of
edges of $G_0$ with $G_1$, say $e_0$ of $G_0$ with $e_1$ of $G_1$, we have
$$W_{(G_1,e_1)}=\gamma\circ W_{(G_0,e_0)}\circ \gamma ^{-1},$$
where we concatenate from right to left, so $\gamma ^2\circ \gamma ^1$ is
defined if and only if $G_1^1=G_0^2$.  For any $\gamma\in \Gamma _I$,
$\{\gamma\circ u_i\circ
\gamma ^{-1}\}_{i=1}^L$ generates $\Isg$ and $\{ \gamma\circ v_j\circ
\gamma ^{-1}\}
_{j=1}^5$ generates $\mathrm{Sp}(2g,\bZ)$.

To complete the proof, suppose that $G_1$ is any trivalent fat graph spine
of
$F$ and $e_1$ is any edge of $G_1$ with distinct endpoints.  There is some
$\varph\in \Msg$ and
$G_0$ with cell in
$D$ so that
$\varph (G_0)=G_1$ since
$D$ is a fundamental domain for $\Msg$, where $e_0$ is the edge of $G_0$
corresponding to
$e_1$.  Writing
$\varph =uv$ (non-canonically), where
$u\in \Isg$ and
$v\in \mathrm{Sp}(2g,\bZ)$, we find that
$$W_{(G_1,e_1)}=u\circ v\circ W_{(G,e)}\circ v^{-1}\circ u^{-1}\in B,$$
from which it directly follows that indeed $B(D)=\mathcal{T}_{g,*}$.
\end{proof}

  From previous computer work \cite{MP} in genus three, there are
12,594 + 42 + 5 generators and
9,548 + 31,760 pentagon and commutativity relations, respectively, plus the
remaining
12,594 involutivity relations.

Taking the generators of 1-3) and adding to the
relations in Theorem \ref{th:2} also the (finitely many) relations of
$\mathrm{Sp}(2g,\bZ)$ gives a finite presentation of
$\Gamma _T$ as one shows again using standard covering space theory.  This
both sheds light
on our presentation of $\Gamma _I$ in
Theorem \ref{th:2} and gives a representation $\mathrm{Sp}(2g,\bZ)\to\Gamma
_T$.

Consider the ``level $N$ classical Torelli group'' ${\mathcal I}_{g,*}(1)(N)$ defined to be
the subgroup of the mapping class group ${\mathcal M}_{g,*}$
that acts identically on the homology of the surface $F$ with coefficients in ${\mathbb Z}/N$.
It follows from the previous theorem that there is also a finite presentation of the
fundamental path groupoid of the ``level $N$ Torelli space'' ${\mathbf T}_1(N)={\mathcal
T}_{g,*}/{\mathcal I}_{g,*}(1)(N)$ with its induced ideal cell decomposition.  
Since the moduli space of $F$ has only finitely many cells (albeit super-exponentially many in
the genus of $F$) and there are only finitely many ${\mathbb Z}/N$ homology markings on a
fatgraph with finitely many edges, it follows that the there are only finitely many cells in the
decomposition of 
${\mathbf T}_1(N)$.  Thus, the group is necessarily of finite index in the groupoid, giving the
following result:

\vskip .2in

\begin{corollary}\label{levelN}
The generators and relations of Theorem~\ref{groupoid} descend to a finite
presentation of the level $N$ Torelli group ${\mathcal I}_{g,*}(1)(N)$, for each $N\geq 1$ and
$g\geq 3$.
\end{corollary}

\vskip .2in

Arguing as in Theorem~\ref{groupoid} and using the  
result of Johnson \cite{Johnson85} 
that the normal subgroup ${\mathcal K}_{g,*}={\mathcal I}_{g,*}(2)$ of ${\mathcal M}_{g,*}$
is normally generated by Dehn twists on $[g/2]$ separating curves, we find the following:

\vskip .2in

\begin{corollary}\label{secondtor}
A generating set for the groupoid corresponding to the Johnson 
subgroup ${\mathcal I}_{g,*}(2)={\mathcal K}_{g,*}$, for $g\geq 3$, is given by adding to 
those in Theorem~\ref{groupoid} the following $[g/2]$ elements:
for each $i=1,\cdots,[g/2]$, choose a Dehn twist on a separating
simple closed curve which cuts the surface $F$ into two pieces
of genera $i$ and $g-i$ and then choose a representing sequence of 
Whitehead moves starting from $G$ corresponding to this element.
Relations in this groupoid are given as in Theorem~\ref{groupoid}.
\end{corollary}

\vskip .2in

\section{Combinatorial realization of the extended Johnson homomorphism}

It was shown in \cite{Morita93} that the Johnson homomorphism
$\tau_1:\Isg\lra \La^3 H$
can be extended to a crossed homomorphism
$\tilde k:\Msg\lra \frac{1}{2} \La^3 H$
essentially uniquely.
Here we prove that it in fact lifts to a {\it canonical} $1$-cocycle $j$
of the cell complex $\hat{\mathcal G}_T$ with values in $\La^3 H$
which is $\Msg$-equivariant. We show that the above
extension of $\tau_1$ to a crossed homomorphism on $\Msg$ is
a simple consequence of the existence of $j$.

\begin{definition}
We define a combinatorial cochain
$$
j\in C^1(\hat{\mathcal G}_T,\La^3 H_1)
$$
by setting
$$
j(W_e)=a\land b\land c=c\land d\land a.
$$
in the notation of Figure~1.
\end{definition}

\begin{theorem}
The cochain $j$ is a cocycle, so that
$
j\in Z^1(\hat{\mathcal G}_T,\La^3 H)
$
and it is $\Msg$-equivariant in the sense that the equality
$$
j(\varph(\sigma))=\varph (j(\sigma))
$$
holds for any oriented $1$-cell $\sigma$ of $\hat{\mathcal G}_T$ and
$\varph\in\Msg$. It follows that $j$ descends to a
$1$-cocycle
$
j\in Z^1(\hat{\mathcal G}_I,\La^3 H)
$
of the cell complex $\hat{\mathcal G}_I$ which is a
$K(\Isg,1)$. Furthermore, the associated
group homomorphism
$$
[j]\in H^1(\hat{\mathcal G}_I,\La^3 H))\cong
\Hom(\Isg,\La^3 H)
$$
coincides with six times the Johnson homomorphism
$\tau_1:\Isg\ra\La^3 H$, i.e.,
$
[j]=6 \tau_1.
$
\label{th:j}
\end{theorem}

\begin{proof}
We begin by checking that $j(W_e)\in\La^3 H$ is a
well-defined $1$-cochain.
For this, we must first check that
$j(W_e)$ is independent of the orientation of the edge
and that involutivity leads to a vanishing sum in $\La^3 H$.
Adopting the notation of Figure~1, we have
$a+b+c+d=0$, and so $j(W_e)=a\land b\land c=c\land d\land a$
is well-defined independent of the orientation.
Likewise for involutivity, we have
$j(W_e)+j(W_f)=a\land b\land c+b\land c\land d=0$.

Next we show that $j$ is a $1$-cocycle.
We have to check that the commutativity and pentagon
relations imply vanishing sums in $\La^3 H$.
Commutivity follows immediately because addition is
commutative in $\La^3 H$.

For the pentagon relation and in the notation of Figure~2,
we have $a+b+c+d+e=0$, and
$j(W_f)+j(W_{g_1})+j(W_{f_2})+j(W_{g_3})+j(W_{f_4})$
is given by
$$
b\wedge c\wedge d ~+~e\wedge a\wedge b ~+~c\wedge d\wedge e
~+~a\wedge b\wedge c~+~d\wedge
e\wedge a =0,
$$
as required.


\vskip .2in

\hskip .6in{\epsffile{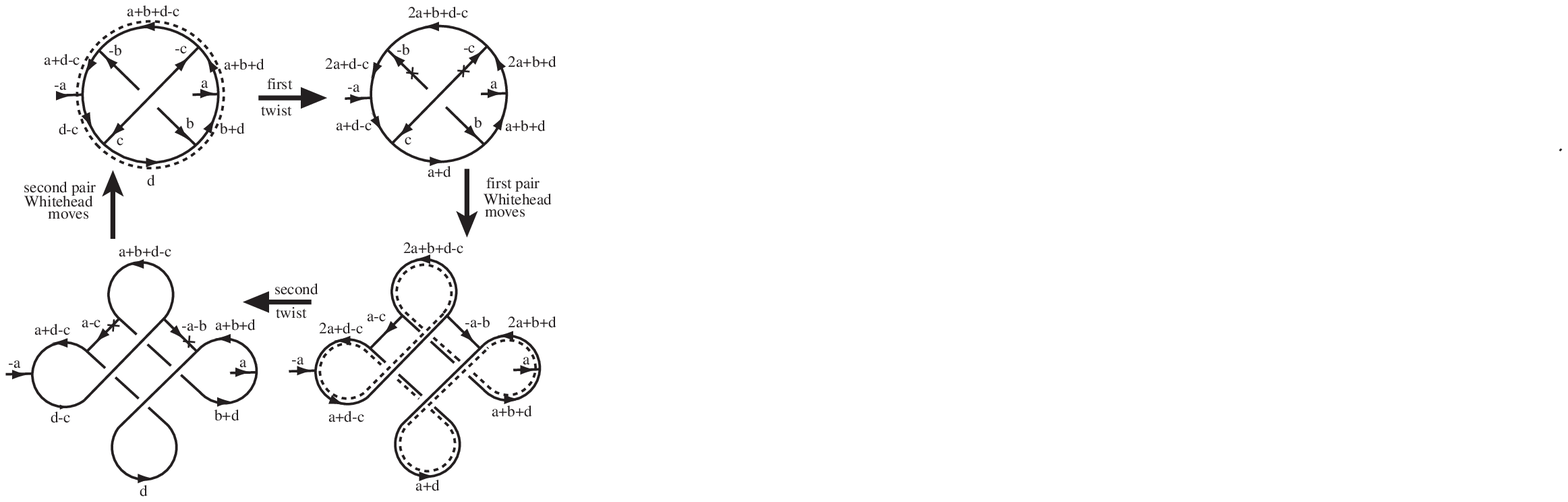}}

\vskip .2in

\centerline{{\bf Figure 3}~The torus $BP$ map $\varph\in\Isg$.}

\vskip .2in

Finally, by \cite{Johnson80}\cite{Johnson85b}, there is a unique homomorphism
$\Isg\ra \La^3 H$ up to multiplication,
and it remains only to find the multiple,
which we do by example:
We compute $j(\varph)$ for a particular $\varph\in\Isg$.

In fact in the terminology of \cite{Johnson83}, $\varph:F\to F$ is a
``torus bounding pair'' mapping, i.e., $\varph$ is defined by Dehn twists
in opposite
directions along the boundary components of a torus-minus-two-disks
embedded in $F$.  A
fat graph spine of this torus-minus-two-disks is drawn in the upper-left of
Figure~3.  There
are four stages to the definition of $\varph$: first, there is a Dehn twist
along the dotted
line in the upper-left to produce the fat graph in the upper-right; second,
there are a pair
of Whitehead moves along the edges marked with crosses in the upper-right
to produce the
fat graph in the lower-right; third, there is another Dehn twist along the
dotted line in the
lower-right to produce the fat graph in the lower-left; fourth and finally,
another pair of
Whitehead moves along the edges marked with crosses in the lower-left to
produce the
fat graph in the upper-left again.  The sequence of Whitehead moves
corresponding to the
first Dehn twist is described on the left hand side of Figure~4,
and the sequence for the second Dehn twist is
described on the right hand side of Figure~4.
In all, the homeomorphism $\varph$ is given by a composition of
$14$ Whitehead moves.


\vskip .2in

\centerline{\epsffile{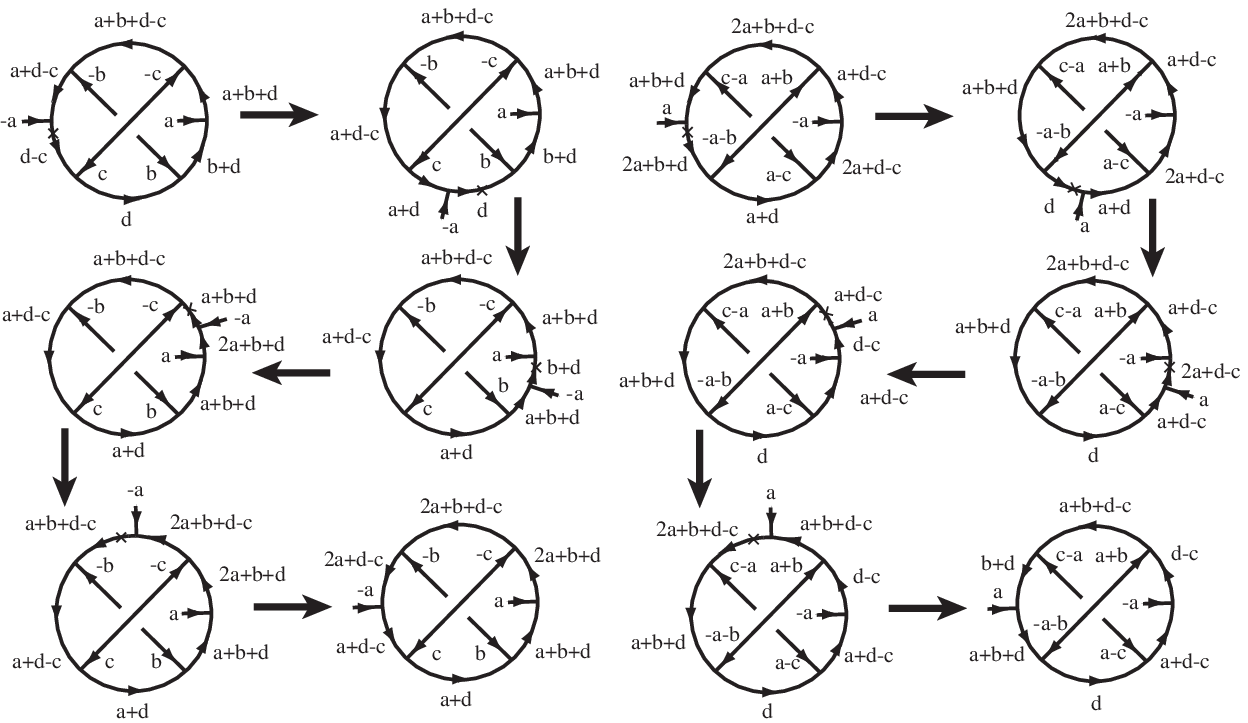}}

\vskip .2in

\hskip .1in{4a}~The first Dehn twist\hskip .7in{4b}~The second Dehn twist

\vskip .2in

\centerline{{\bf Figure 4}~The Dehn twists as Whitehead moves.}

\vskip .2in

To finally compute the explicit value of the crossed homomorphism, we
calculate the contribution to the $3$-form from the first Dehn twist on
the left hand side of Figure~4 to be
\begin{align*}
&-d\wedge c\wedge (a+d-c)~-~(b+d)\wedge b\wedge (a+d)\\
&-(a+b+d)\wedge a\wedge (a+b+d)
+~(a+b+d-c)\wedge c\wedge (2a+b+d)\\
&+~(a+d-c)\wedge b\wedge (2a+b+d-c),
\end{align*}
the contribution from the second Dehn twist on the right hand side in Figure~4 to be
\begin{align*}
&-(a+b+d)\wedge a \wedge (a+d)~-~(2a+d-c)\wedge (a-c)\wedge d\\
&+~(a+d-c)\wedge a\wedge (a+d-c) ~-~(2a+b+d-c)\wedge(a+b)\wedge(d-c)\\
&-~(a+b+d)\wedge (c-a)\wedge (a+b+d-c),
\end{align*}
and the respective contributions from the two pairs of Whitehead moves in
Figure~3 to be
\begin{align*}
&~-~(2a+b+d)\wedge(2a+b+d-c)\wedge (a+d-c)\\
&~-~(a+d)\wedge (a+b+d)\wedge (2a+b+d-c)\\
&~+~(a+b+d-c)\wedge (d-c)\wedge d~+~(a+d-c)\wedge d\wedge (b+d).
\end{align*}
By direct calculation, each of these expressions is equal to
$2~a\wedge b\wedge c$.
Since the value of Johnson's homomorphism on $\varph\in \Isg$ is known
\cite{Johnson80} to be $a\wedge b\wedge c$,
Theorem \ref{th:j} follows for the cochain $j$ on
$\hat{\mathcal G}_T$ with values in $\Lambda^3 H$.
This finishes the proof.
\end{proof}

Now fix a vertex $v_0$ of the complex
$\hat{\mathcal G}_T$. For any element $\varph\in \Msg$,
choose a chain $\sigma_1,\cdots,\sigma_p$ of oriented $1$-cells
which connects $v_0$ to $\varph(v_0)$ and set
$$
\tilde j(\varph)=\sum_{i=1}^p j(\sigma_i).
$$
Since $j$ is a cocyle, this value $\tilde j(\varph)$ does not depend
on the choice of the chain connecting $v_0$ to $\varph(v_0)$, and 
hence we obtain a mapping
$$
\tilde j:\Msg\lra \La^3 H.
$$

\begin{theorem}
The mapping $\tilde j:\Msg\ra \La^3 H$ is a crossed
homomorphism whose restriction to the Torelli group
is six times the Johnson homomorphism.
\label{th:crossed}
\end{theorem}

\begin{proof}
Let $\varph,\psi\in \Msg$ be any two elements.
Choose a chain $\sigma_1,\cdots,\sigma_p$
   (resp. $\tau_1,\cdots,\tau_q$) of oriented
$1$-cells of $\hat{\mathcal G}_T$ which connects
$v_0$ to $\varph(v_0)$ (resp. to $\psi(v_0)$).
Thus, $\varph(\tau_1),\cdots,\varph(\tau_q)$ is a chain which connects
$\varph(v_0)$ to $\varph(\psi(v_0))$, and hence
the composed chain
$\sigma_1,\cdots,\sigma_p, \varph(\tau_1),\cdots,\varph(\tau_q)$ connects
$v_0$ to $\varph\psi(v_0)$.  Furthermore, we have
\begin{align*}
\tilde j(\varph\psi)&=\sum_{i=1}^p j(\sigma_i)+\varph\left(
\sum_{k=1}^q j(\tau_k)\right)\\
&=\tilde j(\varph)+\varph(\tilde j(\psi)),
\end{align*}
using the fact that $j$ is $\Msg$-equivariant.
Thus, $\tilde j$ is a crossed homomorphism as required.
If $\varph$ belongs to the Torelli group $\Isg$, then
the projection of the vertex $\varph(v_0)$ to $\hat{\mathcal G}_I$
is equal to that of $v_0$, and the chain $\sigma_1,\cdots,\sigma_p$
projects to a closed curve on $\hat{\mathcal G}_I$ which represents
the element $\varph\in \Isg=\pi_1 \hat{\mathcal G}_I$.
Since $j$ gives rise to a homomorphism $\Isg\ra\La^3 H$
which coincides with $6 \tau_1$ by Theorem \ref{th:j},
we can conclude the last
claim of the theorem.
\end{proof}

In \cite{Morita93}, there was constructed
a crossed homomorphism
$$
\tilde k: \Msg\lra \frac{1}{2}\La^3 H
$$
whose restriction to the Torelli group
coincides with the Johnson
homomorphism $\tau_1:\Isg\ra \La^3 H$.

\begin{corollary}
We have the identity
$$
[\tilde j]=6 [\tilde k]\in H^1(\Msg;\La^3 H).
$$
\label{cor:crossed}
\end{corollary}

\begin{proof}
This follows from the result proved in \cite{Morita93}
that $H^1(\Msg;\La^3 H)$ is isomorphic to $\bZ^2$, and
furthermore any element in this group is determined
by its restriction to the Torelli group $\Isg$, where it gives
a homomorphism.
\end{proof}

\section{New combinatorial cycles on the moduli space}

In this section by making use of the $1$-cocycle $j$ defined in the
previous section, we construct various combinatorial
cocycles on the cell complex $\hat{\mathcal G}_M$,
which is a spine of the moduli space $\mathbf{M}_{g,*}$ and therefore homotopy
equivalent to $\mathbf{M}_{g,*}$, as well
as various closed cycles on $\mathbf{M}_{g,*}$ which arise as their 
Poincar\'e duals.
We can relate this construction to an earlier
one \cite{Morita96} which was given
in the context of group cohomology of the mapping class group.
By virtue of results of Kawazumi and Morita \cite{KM96}\cite{KM01},
we can conclude that the totality of cohomology classes represented
by the above combinatorial cocycles is precisely the
tautological algebra of the moduli space.

We emphasize that our combinatorial cycles on the moduli space
are completely different from those defined by Witten
and Kontsevich \cite{Kontsevich92},
which were studied first by Penner \cite{Penner93} and then
by Arbarello and Cornalba \cite{AC} and more recently
by Igusa \cite{Igusa1}\cite{Igusa2},
Igusa and Kleber \cite{IK},
and Mondello \cite{Mondello}. The cycles of Witten
and Kontsevich are defined by specifying valencies
of fat graphs embedded in Riemann surfaces, while as we shall see, our
new cycles are defined in terms of homology intersection properties
of simple closed curves dual to edges of fat graphs.

Consider the $1$-cocycle
$
j\in Z^1(\hat{\mathcal G}_T,\La^3 H)
$
and its powers
$$
j^{2k}\in Z^{2k}(\hat{\mathcal G}_T,\La^{2k}\La^3 H)\quad (k=1,2,\cdots).
$$
Since $\hat{\mathcal G}_T$ is not a simplicial complex
but only a cell complex, we must choose a diagonal
approximation to obtain explicit cocycles for the
cup product powers $j^{2k}$. However, any choice
which is $\Msg$-equivariant suffices for the following
argument.
If we apply any $Sp$-invariant contraction
$
C:\La^{2k}\La^3 H\lra \bZ
$
to $j^{2k}$, we obtain a $2k$-cocycle
$
C(j^{2k})\in Z^{2k}(\hat{\mathcal G}_T,\bZ)
$
which is $\Msg$-invariant and hence descends to a
cocycle
$$
C(j^{2k})\in Z^{2k}(\hat{\mathcal G}_M,\bZ)
$$
of the cell complex $\hat{\mathcal G}_M,$.
Taking the Poincar\'e dual of such a
cocycle, we obtain a closed cycle of $\mathbf{M}_{g,*}$
with respect to its canonical cell decomposition.

Recall from \cite{Morita96} that we can enumerate all
$Sp$-invariant contractions $\La^*\La^3 H\ra \bZ$ as
follows. Let ${\mathcal G}$ denote the set of all the isomorphism
classes of trivalent graphs $\vGa$. There is a surjective
homomorphism
$$
\bQ[\vGa;\vGa\in {\mathcal G}]\lra \Hom(\La^*\La^3 H,\bQ)^{Sp},
$$
where the left hand side denotes the polynomial algebra generated
by ${\mathcal G}$. We denote by $C_\vGa$ (resp. $Z_\vGa$)
the cocycle of the cell complex $\hat{\mathcal G}_M$
(resp. cycle of the moduli space) corresponding to
$\vGa\in{\mathcal G}$.
Thus, we obtain the following result.

\begin{theorem}
For each trivalent graph $\vGa$ with $2k$ vertices,
there corresponds a $2k$-cocycle $C_\vGa$ of the
cell complex $\hat{\mathcal G}_M$, which induces a
homomorphism
$$
\Phi: \bQ[\vGa;\vGa\in {\mathcal G}]\lra
H^*(\hat{\mathcal G}_M;\bQ)\cong H^*({\mathbf M}_{g,*};\bQ).
$$
\end{theorem}

The homomorphism $\Phi$ in the above theorem can be considered
as the combinatorial realization of the homomorphism
$$
\Phi': \bQ[\vGa;\vGa\in {\mathcal G}]\lra
H^*(\Msg;\bQ)
$$
defined in \cite{Morita96} (see Theorem 3.4).
More precisely, we have the following result.

\begin{theorem}
The homomorphism $\Phi$ coincides with the constant $6^*$
times $\Phi'$
through the canonical isomorphism
$H^*({\mathbf M}_{g,*};\bQ)\cong H^*(\Msg;\bQ)$.
\label{th:phiphi}
\end{theorem}

Before proving this theorem, we must prepare a technical
result which should be considered as a standard fact about
the cohomology groups with values in a local coefficient system
(see \cite{Eilenberg}).

Let $K$ be a cell complex which is an Eilenberg-MacLane space
of type $K(\pi_1 K,1)$.
Assume that there is given
a homomorphism
$$
\rho:\pi_1 K\lra \mathrm{Aut}\, L
$$
where $L$ is an abelian group.
There is then a corresponding local coefficient system ${\mathcal L}$
over $K$, well-defined up to isomorphism,
whose characteristic homomorphism coincides
with $\rho$. If there is given a {\it crossed homomorphism}
\begin{equation}
h:\pi_1 K\lra L,
\label{eq:crossed}
\end{equation}
namely, if the equality $h(\ga\ga')=h(\ga)+\ga h(\ga')$
holds for any $\ga,\ga'\in\pi_1 K$,
then it defines an element
$$
[h]\in H^1(\pi_1 K;L)
$$
in the first cohomology of the group $\pi_1 K$ with
coefficients in the $\pi_1 K$-module $L$.
On the other hand, since $K$ is an Eilenberg-MacLane space
of type $K(\pi_1 K,1)$ by assumption, we have
a canonical isomorphism
\begin{equation}
H^*(\pi_1 K;L)\cong H^*(K;{\mathcal L}).
\label{eq:L}
\end{equation}
On the right hand side of the above isomorphism,
the crossed homomorphism \eqref{eq:crossed}
can be realized combinatorially as follows.

Choose a vertex $p_0\in K$ as a base point of $K$
and also a maximal tree $T$ in the $1$-skeleton of $K$.
Any oriented $1$-cell $\sigma\in K$ then defines an element
$\ga_\sigma\in \pi_1 K$ (relative to the base point $p_0$ as well
as the maximal tree $T$). We define a mapping
$$
h':\{\text{oriented $1$-cell of $K$}\}\lra L
$$
by setting
$$
h'(\sigma)=h(\ga_\sigma)\in L.
$$
Note that $h'(\sigma)=0$ for any oriented $1$-cell contained in $T$
and also that $h'(\bar\sigma)=h(\ga_{\bar\sigma})=-\ga_\sigma^{-1}
h'(\sigma)$, where
$\bar\sigma$ denotes the same $1$-cell $\sigma$
but equipped with the opposite orientation.
It is now easy to see that the above $h'$ is in fact a $1$-cocycle
of $K$ with values in the local coefficient system ${\mathcal L}$,
namely, we have
$h'\in Z^1(K;{\mathcal L})$.

There is another combinatorial realization of the crossed
homomorphism \eqref{eq:crossed} used in this paper, which
can be described as follows.

Let $\widetilde{K}$ be the universal covering cell complex of $K$
and let $\pi:\widetilde{K}\to K$ be the projection.
Since $\widetilde{K}$ is simply connected, the local system
$\pi^*{\mathcal L}$ is trivial, namely, it is isomorphic to
the constant local system $L$ over $\widetilde{K}$.
In fact, if we choose a point $\tilde p_0\in\widetilde{K}$
which projects to the base point $p_0\in K$,
we obtain an explicit isomorphism
$\pi^*{\mathcal L}\cong L$, and it is then easy to see that the
pull back homomorphism
$$
\pi^*:Z^1(K;{\mathcal L})\lra
Z^1(\widetilde{K};\pi^*{\mathcal L})
\cong Z^1(\widetilde{K};L)
$$
induces an isomorphism
$$
Z^1(K;{\mathcal L})\overset{\pi^*}{\cong}
Z^1_{\pi_1 K}(\widetilde{K};L),
$$
where the right hand side denotes the subspace consisting of
all the $1$-cocycles $c\in Z^1(\widetilde{K};L)$ such that
the equality
$$
c(\ga \tilde\sigma)=\ga c(\tilde\sigma)
$$
holds for any $\ga\in\pi_1K$ and any oriented $1$-cell $\tilde\sigma$
of $\widetilde{K}$.
In particular, the pull back $1$-cocycle
$\pi^* h'\in Z^1_{\pi_1 K}(\widetilde{K};L)$
is a $\pi_1 K$-equivariant $1$-cocycle of $\widetilde{K}$.

Conversely, suppose that there is given a $\pi_1 K$-equivariant
$1$-cocycle
$$
\tilde h'\in Z^1_{\pi_1 K}(\widetilde{K};L).
$$
On the one hand,
it projects to a $1$-cochain
$$
h'\in C^1(\widetilde{K};L),
$$
while on the other hand, there corresponds a crossed homomorphism
$$
h:\pi_1 K\lra L
$$
defined by
$$
h(\ga)=\sum_{i=1}^p \tilde h'(\tilde\sigma_i)\quad (\ga\in\pi_1 K),
$$
where $\tilde\sigma_1,\cdots,\tilde\sigma_p$ is a chain of oriented $1$-cells
of $\widetilde{K}$ which connects $\tilde p_0$ to $\ga(\tilde p_0)$.
The fact that $h$ is indeed a crossed homomorphism follows from
the same argument as in the proof of
Theorem \ref{th:crossed}.
We summarize this discussion in the following statement.

\begin{proposition}
In the notation above, the correspondences
\begin{align*}
&\text{crossed homomorphism $h:\pi_1 K\ra L$}\\
\leftrightarrow\ &
\text{$1$-cocycle $h'$ on $K$ with values in the local coefficient system
$\mathcal L$}\\
\leftrightarrow\ &
\text{$\pi_1 K$-equivariant $1$-cocycle $\tilde h'$ of $\widetilde{K}$
with values in $L$}
\end{align*}
realize the isomorphisms
$$
H^1(\pi_1 K;L)\cong H^1(K;{\mathcal L})\cong
H^1(\mathrm{Hom}_{\pi_1 K}(C_*(\widetilde{K}),L))
$$
on the level of cocycles.
\label{prop:local}
\end{proposition}

\noindent
\textit{Proof of Theorem \ref{th:phiphi}}.
We must relate the construction of $\Phi$ to that of $\Phi'$.
The former is defined by taking various contractions of
powers of the $1$-cocycle
$$
j\in Z^1_{\Msg}(\hat {\mathcal G}_T;\La^3 H),
$$
which is an $\Msg$-equivariant $1$-cocycle on the cell
complex $\hat {\mathcal G}_T$.
On the other hand,
the latter is obtained by the same procedure
but instead of $j$, by making use of the crossed homomorphism
$$
\tilde k: \Msg\lra \frac{1}{2}\La^3 H
$$
constructed in \cite{Morita93}.
By Corollary \ref{cor:crossed} and
Proposition \ref{prop:local}, the proof would be complete
if $\mathbf{M}_{g,*}$ were an Eilenberg-MacLane
space of type $K(\Msg,1)$.
This is not the case, however,
because the mapping class groups $\Msg$ have torsion
so that the action of $\Msg$ on $\La^3 H$
does not induce a local coefficient
system on the moduli space $\mathbf{M}_{g,*}$ in the usual
sense but only in an orbifold sense.

However, we can proceed as follows. Choose a torsion free
normal subgroup $\mathcal{N}\subset\Msg$ of finite index,
e.g., the kernel of the natural representation
$\Msg\ra\mathrm{Sp}(2g,\bZ/3)$,
and let $\mathcal{F}$ denote the quotient finite group $\Msg/\mathcal{N}$.
The quotient space $\mathbf{M}_{\mathcal N}=\mathcal{T}_{g,*}/\mathcal{N}$
is an Eilenberg-MacLane space of type $K(\mathcal{N},1)$, and hence
we have an isomorphism
\begin{equation}
H^1(\mathcal{N};\La^3 H)\cong H^1(\mathbf{M}_{\mathcal N};\La^3\mathcal{H}),
\label{eq:n}
\end{equation}
where $\La^3\mathcal{H}$ denotes a local coefficient system
on $\mathbf{M}_{\mathcal N}$ induced by the representation
$\mathcal{N}\subset\Msg\ra\mathrm{Sp}(2g,\bZ)$.
Under the above isomorphism \eqref{eq:n},
the restriction of the class $[\tilde k]$ to the left hand side
corresponds to that of $6 [\tilde j]$ to the right hand side
by Corollary \ref{cor:crossed} and
Proposition \ref{prop:local} suitably adapted to the
context of the subgroup $\mathcal{N}$.
If we take contractions of powers of these classes,
it is clear that the resulting cohomology classes
belong to the $\mathcal{F}$-invariant subspaces
$$
H^*(\mathcal{N};\bQ)^\mathcal{F}\cong H^*(\mathbf{M}_{\mathcal
N};\bQ)^\mathcal{F}.
$$
On the other hand, we have canonical isomorphisms
\begin{align*}
H^*(\Msg;\bQ)&\cong H^*(\mathcal{N};\bQ)^\mathcal{F},\\
H^*(\mathbf{M}_{g,*};\bQ)&\cong H^*(\mathbf{M}_{\mathcal N};\bQ)^\mathcal{F}
\end{align*}
because the finite group $\mathcal{F}=\Msg/\mathcal{N}$ acts on
$\mathbf{M}_{\mathcal N}$
with quotient the moduli space $\mathbf{M}_{g,*}$.
  From this, it finally follows that $\Phi=6^* \Phi'$ as required.
\qed

\vskip .2in

Since it was proved in \cite{KM96} that the image of the homomorphism
$\Phi'$ is the tautological algebra
${\mathcal R}^*(\Msg)$ of $\Msg$, which is by definition
the subalgebra of $H^*(\Msg;\bQ)$ generated by the Mumford-Morita-Miller
classes $e_i$ and the universal Euler class $e\in H^2(\Msg;\bZ)$,
we obtain the following corollary to Theorem \ref{th:phiphi}.

\begin{corollary}
The image of the homomorphism $\Phi$ is precisely the tautological
algebra of the moduli space ${\mathbf M}_{g,*}$.
\end{corollary}

Next, we compute our cocycles in codimension two explicitly.

\begin{example}
The cup square $j^2\in Z^2(\hat{\mathcal G}_T;\La^2\La^3 H)$ can be computed
as follows. There are two kinds of $2$-cells in $\hat{\mathcal G}_T$.
One is the square as shown in Figure~5 and the other is the
pentagon
depicted in Figure~6.


\vskip .2in

\centerline{\epsffile{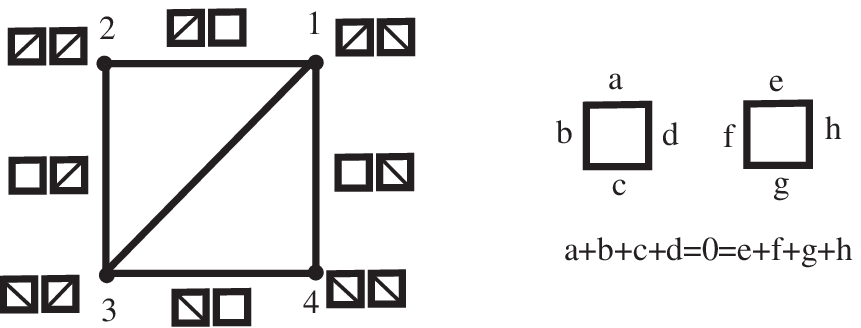}}

\vskip .2in

\centerline{{\bf Figure 5}~The square.}

\vskip .2in

We first consider the simpler case of two complementary quadrilateral and
adopt the notation of Figure~5, where the right hand
side of the figure shows the homology marking, and the left hand
side shows a $2$-cell in $\hat{\mathcal G}_T$ together with a
triangulation of this $2$-cell. Using the Alexander-Whitney diagonal 
approximation, we may compute
\begin{align*}
j^2(1,2,3)+j^2(1,3,4)=& j(2,3)\wedge j(1,2)+j(3,4)\wedge j(1,3)\\
=& j(2,3)\wedge j(1,2) +j(3,4)\wedge \{j(1,4)+j(4,3)\}\\
=& j(2,3)\wedge j(1,2)+j(3,4)\wedge j(1,4)\\
=& abc\wedge fgh + efg\wedge abc\\
=& 2~efg\wedge abc.
\end{align*}
As to invariance, we may triangulate the $2$-cell on the left hand side
in Figure~5 with the other diagonal and likewise calculate
that $j^2(1,2,4)+j^2(2,3,4)=2~efg\wedge abc$ independent of triangulations.


\vskip .2in

\centerline{\epsffile{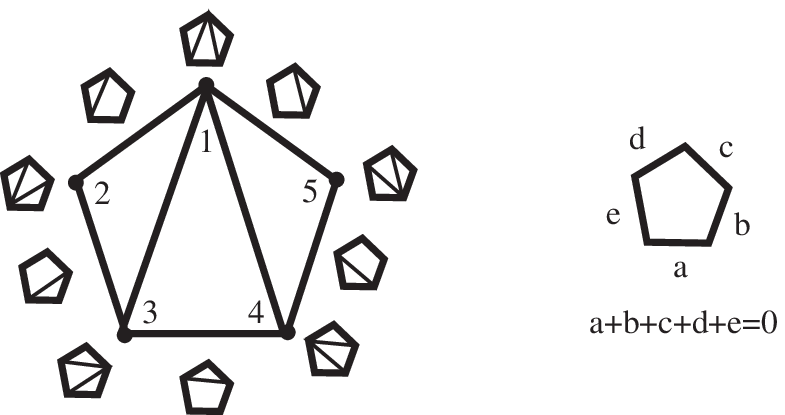}}

\vskip .2in

\centerline{{\bf Figure 6}~The pentagon.}

\vskip .2in

For the pentagon,
adopt the notation of Figure~6, and similarly
calculate
\begin{align*}
j^2=& j^2(1,2,3)+j^2(1,3,4)+j^2(1,4,5)\\
=& j(2,3)\wedge j(1,2)+j(3,4)\wedge j(1,3)+j(4,5)\wedge j(1,4)\\
=& j(2,3)\wedge j(1,2)+j(3,4)\wedge\{j(1,2)+j(2,3)\}
                        +j(4,5)\wedge\{j(1,5)+j(5,4)\}\\
=& j(2,3)\wedge j(1,2)+j(3,4)\wedge\{j(1,2)+j(2,3)\}
                        +j(4,5)\wedge j(1,5).
\end{align*}
  From the figure and the definition of $j$, we may calculate
\begin{align*}
j(1,5)=& dea = -j(5,1)=ea(b+c),\\
j(5,4)=& bcd = -j(4,5)=cd(a+c),\\
j(4,3)=& eab = -j(3,4)=ab(c+d),\\
j(2,1)=& abc = -j(1,2)=be(d+c),\\
j(3,2)=& cde = -j(2,3)=de(a+b),
\end{align*}
and substituting into the expression for $j^2$, we find
$$
j^2=cde\wedge abc + eab\wedge abc + eab\wedge cde + dea\wedge bcd.
$$
We may also eliminate $e=-a-b-c-d$ to get
\begin{align*}
j^2= adb&\wedge abc +dab\wedge cda + bcd\wedge dca\\
    &+2\{dca\wedge abc +cba\wedge abc + dab\wedge cdb\}.
\end{align*}
As to invariance, we may triangulate the pentagonal $2$-cell on the left 
hand side in Figure~6 so
that both diagonals have vertex $5$ as endpoint, and similarly calculate
that $j^2=j^2(1,2,5)+j^2(2,3,5)+j^2(3,4,5)$. Upon again eliminating
$e$, we find the same expression for $j^2$.

Thus, we find that the value of the pull back of $j^2$ on any $2$-cell in 
$\hat{\mathcal G}_T$ is
{\it independent}  of the
triangulation of the $2$-cell in either case.

Consider the two trivalent graphs $\vGa_i\ (i=1,2)$ with two vertices
where $\vGa_1$ has two loops while $\vGa_2$ is the graph isomorphic to the 
letter theta.
We have the corresponding two contractions
$$
C_1,C_2: \La^2\La^3 H\lra \bZ
$$
which are explicitly given by
\begin{align*}
C_1\bigl(&(a_1\land a_2\land a_3)\land (b_1\land b_2\land b_3)
\bigr),\\
=&
\sum_{\sigma,\tau\in {\mathfrak S}_3}
\sgn\sigma \sgn\tau\,
(a_{\sigma(1)}\cdot a_{\sigma(2)})(b_{\tau(1)}\cdot b_{\tau(2)})
(a_{\sigma(3)}\cdot b_{\tau(3)})\\
C_2\bigl(&(a_1\land a_2\land a_3)\land (b_1\land b_2\land b_3)
\bigr),\\
=&
\sum_{\sigma,\tau\in {\mathfrak S}_3}
\sgn\sigma \sgn\tau\,
(a_{\sigma(1)}\cdot b_{\tau(1)})(a_{\sigma(2)}\cdot b_{\tau(2)})
(a_{\sigma(3)}\cdot b_{\tau(3)})
\end{align*}
where $a_i, b_i \in H$ and $\cdot$ denotes the homology intersection pairing.
If we apply these contractions to $j^2$, we obtain two
cocycles $C_{\vGa_i}$ of the cell complex $\hat{\mathcal G}_M$.
The results of \cite{Morita89b}\cite{Morita96} together with
Theorem \ref{th:phiphi} imply
\begin{align*}
\frac{1}{36}[C_{\vGa_1}]&= -4g(g-1)e-e_1\\
\frac{1}{36}[C_{\vGa_2}]&= 6 g e-e_1.
\end{align*}
\end{example}

The analogous calculation of higher powers of $j$ is clear, and
Theorem \ref{th:phiphi} implies that
the cohomology class $\Phi(\vGa)$ is the same as $6^{2k} \Phi'(\vGa)$
for any trivalent graph $\vGa$ with $2k$ vertices.
The explicit formula for the latter
was given in Kawazumi and Morita \cite{KM01}.

\section{Nilpotent markings and presentations of higher Torelli groups}

Let $F$ denote a closed oriented surface with base point $*\in F$
and let $F'$ denote the punctured surface $F-\{*\}$ as before.

\begin{definition}
[$N_k$- and $\pi _1$-markings]
Given a trivalent fat graph spine $G$ of
$F'$ and some $k\geq 1$, define an $N_k$-{\it marking} on $G$ to be an
assignment
$a\mapsto \alpha _a\in N_k$ to each oriented edge $a$ of $G$ so that:

\begin{enumerate}
\item
if $\alpha$ is assigned to the oriented edge $a$ of $G$, then
$\bar\alpha=\alpha ^{-1}$
is assigned to the reverse of $a$;
\item
given three oriented edges $a,b,c$ of $G$ pointing towards a
common vertex $v$ in this counter-clockwise cyclic order, we have
$$
\alpha _a ~\alpha _b ~\alpha _c =1\in N_k;
$$
\item
the marking is of ``full rank'' in the sense that
$$
\{ \alpha
_a: a~{\rm is~an~oriented~edge~of}~G\}
$$
generates $N_k$.
\end{enumerate}
\end{definition}

Likewise, a $\pi _1$-{\it marking} on $G$ is an assignment
$a\mapsto \alpha_a\in\pi_1(F,*)$ satisfying the analogous properties 1-3).

Thus, a $\pi _1$-marking on $G$ determines
a residual $N_k$-marking for each $k\geq 1$, and furthermore, an
$N_\ell$-marking on $G$
determines a residual $N_k$-marking provided $\ell\geq k$.

\begin{definition}[tautological $\pi_1$-marking]
Let $G$ be any trivalent fat graph spine of $F'$
so that it defines a $(6g-4)$-cell of the Teichm\"uller space
${\mathcal T}_{g,*}$ and also its Poincar\'e dual is a vertex of
$\hat{\mathcal G}_T$.
An oriented edge $a$ of $G$ with dual arc
$a'\in A_G$ determines an element
$\alpha _a\in\pi _1(F,*)$ by the following procedure:
Orient $a'$ so that the
orientation given by $a',a$ in this order agrees with the orientation of
$F$ at $a\cap a'$.
The oriented curve $a'$ connects the base point $*$
to itself and  determines an element in $\pi _1(F,*)$.
This assignment
$a\mapsto \alpha _a\in\pi_1(F,*)$
is said to be {\it induced} by $G$ in $F$.
\end{definition}

\begin{lemma}
Given a trivalent fat graph spine $G$ of $F'$, the induced assignment
$a\mapsto \alpha _a\in \pi _1(F,*)$ is a $\pi _1$-marking on $G$.
\label{lem:pi}
\end{lemma}

\begin{proof}
Conditions 1) and 2) follow by construction, and 3) follows since $G$ is a
spine of $F'$.
\end{proof}

We call the above $\pi_1$-marking the {\it tautological marking}.
This marking induces an $N_k$-marking for any $k\geq 1$ which we also
call the {\it tautological} $N_k$-marking.

Now we have the following result.

\begin{theorem}
The higher Torelli group $\Isg(k)$, for $k\geq 1$, admits
the following infinite presentation.  Generators are given by sequences of
Whitehead
moves starting from a fixed trivalent fat graph spine $G$ of $F'$ with its
tautological $N_k$-marking and
ending with $\varph (G)$, for
$\varph\in \Msg$, which leave invariant the $N_k$-marking.  Relations are
given by
involutivity and the commutativity or pentagon relations corresponding to
the codimension-two
faces in ${\mathcal G}_T/\Isg(k)$.
\end{theorem}

\begin{proof}
Since cells in ${\mathcal G}_T/\Isg(k)$ are indexed by $\Isg(k)$-orbits of
$N_k$-marked fat graphs, the proof is analogous to that of
Theorem \ref{th:presentation}.
\end{proof}

\vskip .2in

We shall finally discuss the Johnson homomorphisms
$$\tau _k:\Isg(k)\to \mathrm{Hom}(N_{k+1},\Ga_k/\Ga_{k+1})
$$
in this context for $k\geq 1$.

Suppose there is a sequence of Whitehead moves
starting from a trivalent fat graph
$G$ with its tautological $N_{k+1}$-marking $\mu$ and ending with the trivalent
fat graph $G'$ and its
corresponding $N_{k+1}$-marking $\mu '$, and let $e\mapsto e'$ denote the 
correspondence of edges of
$G$ and $G'$ induced by the Whitehead moves.

Suppose there is an element $\varph\in\Isg(k)$ so that
$\varph (G)=G'$ and $\varph (e)=e'$ for each edge $e$ of $G$, i.e., suppose 
that the sequence of
Whitehead moves corresponds to some $\varph\in\Isg(k)$.  In this case, we 
may define the element
$$\lambda _k(e)=\mu (e) \bigl [\mu '(e' )\bigr]^{-1}~\in ~N_{k+1},~{\rm 
for}~e~{\rm an~edge~of}~G.$$
By definition of $\Isg(k)$, it follows that $\lambda _k(e)$ maps to zero in $N_k$, i.e., 
$\lambda _k(e)\in\Ga_k$ for
each edge $e$, and according to the basic exact sequence (1), we may thus 
regard $$\lambda
_k(e)\in\Ga_k/\Ga_{k+1}~{\rm for~each~edge}~e.$$

\noindent Suppose that there is a trivalent vertex of $G$ with incident 
edges $a,b,c$ in this correct
counter-clockwise order, where each edge is oriented pointing towards their 
shared vertex.  Abusing
notation slightly letting $x=\mu (x)$, $x'=\mu ' (x')$ and setting 
$\overline x =x^{-1}$,
for $x=a,b,c$, there
is a small but satisfying calculation in $N_{k+1}$ to see that
$\lambda _k$ determines a homomorphism:

\begin{align*}
(a\overline a')(b\overline b')(c\overline c')=&(a(b\overline b')\overline 
a')(c\overline c') \\
=&(a(b(c\overline c')\overline b')\overline a')\\
=&(abc)(\overline c'\/\overline b'\/\overline a')\\
=&1,\\
\end{align*}

\noindent
since $abc=1=a'b'c'$ in $N_{k+1}$ by definition of $N_{k+1}$-marking and 
using the fact
that elements of $\Gamma_k$ such as $x\bar x'$ commute in in $N_{k+1}$ .
The next
result formalizes this discussion since $\lambda _k$ agrees with $\tau _k$ 
on generators by
definition.

\begin{proposition}
The combinatorially defined mapping $$\lambda_k:\Isg(k)\to 
\mathrm{Hom}(N_{k+1},\Ga_k/\Ga_{k+1})$$
in the notation above agrees with the Johnson homomorphism $\tau _k$ for 
each $k\geq 1$.
\end{proposition}

\section{Concluding remarks}

It was proved in \cite{Morita96} that the natural action of
$\Msg$ on the third nilpotent quotient $N_3$
of $\pi_1(F,*)$ can be explicitly
described as a homomorphism
\begin{equation}
\rho_3:\Msg\lra \left(\frac{1}{24} H_2\tilde{\times}
\frac{1}{2}\Lambda^3 H\right)\rtimes \mathrm{Sp}(2g,\bZ).
\label{eq:rho3}
\end{equation}
Here $H_2$ denotes a certain space of
$4$-tensors of $H$, $\tilde{\times}$ denotes the central
extension of $\frac{1}{2}\Lambda^3 H$ by
$\frac{1}{24} H_2$ given by a bilinear
skew symmetric pairing
$$
\Lambda^3 H\otimes \Lambda^3 H\lra
H_2
$$
(see \cite{Morita96} for details).

Ignoring the integral denominators, there is a forgetful homomorphism
$p:H_2\tilde{\times}\Lambda^3 H\to \Lambda^3 H$ and this induces
a {\it mapping}
$$
H^1(\Msg;H_2\tilde{\times}\Lambda^3 H)\overset{p_*}{\lra} 
H^1(\Msg;\Lambda^3 H),
$$
where the left hand side denotes the first cohomology {\it set}
of the group $\Msg$ with values in the non-abelian group
$H_2\tilde{\times}\Lambda^3 H$. As was already mentioned in $\S3$,
there is a canonical element $[\tilde k]\in H^1(\Msg;\Lambda^3 H)$
which extends the first Johnson homomorphism $\tau_1$.

The existence of the representation \eqref{eq:rho3},
or more precisely the fact that the image of the
representation $\Msg\to \mathrm{Aut}(N_3)$ is a {\it split
extension} over $\bQ$, implies that there exists $k_2\in 
H^1(\Msg;H_2\tilde{\times}\Lambda^3 H)$
such that $p_*(k_2)=[\tilde k]$.

One of the principal results of this paper is that there exists a completely
canonical $1$-cocycle, namely
$j\in Z^1_{\Msg}(\hat{\mathcal G}_T;\Lambda^3
H)$,
which realizes the element $[\tilde k]$ combinatorially.
Note here that there is no canonical choice of a crossed homomorphism
$\tilde k$
at the cocycle level, and the uniqueness holds only after taking the
cohomology class, while our $j$ is already canonical at the cocycle level.
In view of this result, it seems reasonable to expect that there
should exist a {\it canonical} element
$$
j_2\in Z^1_{\Msg}(\hat{\mathcal G}_T;H_2\tilde{\times}\Lambda^3 H)
$$
which projects to $j$ under the obvious forgetful mapping.
Note here that again the existence of \eqref{eq:rho3} itself
implies the existence of such an element, say by choosing a maximal tree
in a fundamental domain of $\mathcal{T}_{g,*}$ with respect to the
action of $\Msg$ on it. However, the problem here is whether we can
construct such an element independently of any choices.

At the level of the moduli space $\mathbf{M}_{g,*}$, the above
problem could be phrased as asking for the existence of a certain
canonical element in the left hand side of
$$
Z^1(\mathbf{M}_{g,*};\mathcal{H}_2\tilde{\times}\Lambda^3 \mathcal{H})
\lra Z^1(\mathbf{M}_{g,*};\Lambda^3 \mathcal{H})
$$
which projects to $j$, where $\mathcal{H}_2\tilde{\times}\Lambda^3 \mathcal{H}$
denotes the non-abelian local system over the moduli space
(in the obifold sense) corresponding to
the group $H_2\tilde{\times}\Lambda^3 H$ on which $\Msg$ acts naturally.

We can go further to consider the action of $\Msg$ on higher nilpotent
quotients $N_k$ of $\pi_1(F,*)$ and ask for the existence of a series of
{\it canonical} elements
$$
j_k\in Z^1_{\Msg}(\hat{\mathcal G}_T;M_k)
$$
where $M_k$ denotes a certain nilpotent group sitting inside
$\mathrm{Aut}(N_{k+1})$.
One result in the fundamental work of Hain \cite{Hain} is
that the image of $\Msg$ in $\mathrm{Aut}{N_k}$ is
a split extension over $\bQ$ for any $k$; this implies the existence of such
elements, and the problem here is again whether they can be canonically 
constructed.
In view of Hain's further results in the above cited
paper, and also Kawazumi's recent work in \cite{Kawazumi05},
it seems reasonable to expect an affirmative solution to this problem.
We shall discuss these problems further in a forthcoming paper.

In analogy to the last section of the paper, one can consider the surface
$F'=F-\{*\}$ {\it with base point} $p\in F'$, and ask for presentations of 
the Torelli groups of
these punctured surfaces with base point.  In order to implement the
representations of the higher Torelli groups in this context, one must 
``solve the base point
problem'' and keep track of the surface base point $p$ in the punctured 
surface $F'$ during the
combinatorial moves.  In fact, one can
introduce ``enhanced Whitehead moves'' which act not only on $N_k$-markings 
but also
simultaneously on Penner's coordinates on Teichm\"uller space from 
\cite{Penner87} in order to
describe the evolution of the surface base point $p$ under these enhanced 
moves; this renders these
presentations of higher Torelli groups for punctured surfaces explicitly 
calculable.  This too will
be discussed in a forthcoming paper.

\bibliographystyle{amsplain}

\end{document}